\documentclass[12pt]{article}
\usepackage{amsmath,amssymb,amsthm}
\usepackage[top=20truemm,bottom=25truemm,left=20truemm,right=20truemm]{geometry}
\newtheorem{theorem}{Theorem}[section]

\newtheorem{lemma}[theorem]{Lemma}
\newtheorem{remark}[theorem]{Remark}

\allowdisplaybreaks[1]
\begin{document}

	\title{Certain weighted sum formulas for multiple zeta values 
			with some parameters}
	\author{Shin-ya Kadota}
	\date{}
	\maketitle


	\begin{abstract}
		Multiple zeta values (MZVs) are real numbers which are defined by certain multiple series.
		Recently, many people have researched for relations among them and many relations are well known.
		In this paper, we get a new relation among them 
		which is a generalization of a formula obtained by Eie, Liaw and Ong \cite[Main Theorem]{ELO}, and has four parameters.
		Moreover, by using similar method, we also obtain a relation which has six parameters.
		
	\end{abstract}

\section{Introduction and the statement of main results}

			For an $n$-tuple of  natural numbers ${\boldsymbol\alpha}=(\alpha_1,\alpha_2,\cdots,\alpha_n)~(\alpha_n\geq2)$,
		 	 the multiple zeta value (MZV) $\zeta(\boldsymbol\alpha)=\zeta(\alpha_1,\alpha_2,\cdots,\alpha_n)$ is defined as the following series.
				\[
					\zeta(\boldsymbol\alpha)=\zeta(\alpha_1,\alpha_2,\cdots,\alpha_n):=
					\mathop{\sum}\limits_{0<{\ell}_1<{\ell}_2<\cdots<{\ell}_n}
					\frac{1}{{{\ell}_1}^{\alpha_1}{{\ell}_2}^{\alpha_2}\cdots{{\ell}_n}^{\alpha_n}}.
				\]

	We call $n$ the depth of $\zeta(\boldsymbol\alpha)$, $m=\alpha_1+\cdots+\alpha_n$
	the weight of $\zeta(\boldsymbol\alpha)$.
	Easily one find that MZV has the following integral representation:
		\begin{align*}
			\zeta(\alpha_1,\alpha_2,\cdots,\alpha_n)
			=\mathop{\int\cdots\int}\limits_{0<t_1<t_2<\cdots<t_{m}<1}
				\omega_{1}(t_1)\omega_{2}(t_2)\cdots\omega_{m}(t_m),
		\end{align*}
	where
		\[
			\omega_i(t)=
				\begin{cases}
					\displaystyle\frac{dt}{1-t}&,~i\in\{1,\alpha_1+1,\alpha_1+\alpha_2+1,\cdots,\alpha_1+\cdots+\alpha_{n-1}+1\},\\
					\displaystyle\frac{dt}{t}&,~{\rm otherwise}.
				\end{cases}
		\]

	This integral representation plays an important role in this paper.\par

	First, MZVs were studied by Euler \cite{E} in the special case $n=2$, and Hoffman \cite{H} gave the above general definition.
	It is known that MZVs are closely related to knot theory, arithmetic geometry and mathematical physics.  
	Zagier proposed the conjecture for the dimension of the $\mathbb{Q}$--vector spaces ${\cal Z}_k$ spanned by MZVs with weight $k$ by numerical calculations.
	According to that conjecture we can expect that there are many relations among MZVs.
	So in the study of MZVs, one of the main topic is to obtain various relations among them.
	The present paper is devoted to this topic, and the aim of the present paper is to prove the following two relations.

	Let ${\mathfrak S}_n$ be symmetric group of n-th order.
		\begin{theorem}
				\label{theorem:sub theorem}
		For non-negative integers $k$ and $\ell$, and four parameters (indeterminates) $\mu_1,\mu_2,\xi_1$ and $\xi_2$,
		we have 
			\begin{align*}
				&\sum_{
					\begin{subarray}{c}
						a_1+a_2=k\\
						b_1+b_2=\ell
					\end{subarray}}
					\mu_1^{a_1}\mu_2^{a_2}\xi_1^{b_1}\xi_2^{b_2}
					\zeta(a_1+b_1+2)\zeta(a_2+b_2+2)\\
				=&\sum_{\sigma\in {\mathfrak S}_2}\bigg[
				\sum_{
					\begin{subarray}{c}
						a_1+a_2=k\\
						b_1+b_2=\ell
					\end{subarray}}
				\mu_{\sigma(1)}^{a_1}\mu_{\sigma(2)}^{a_2}
				\xi_{\sigma(1)}^{b_1}\xi_{\sigma(2)}^{b_2}
				\sum_{
					\begin{subarray}{c}
						|\boldsymbol\alpha|=a_1+b_1+1\\
						|\boldsymbol\beta|=a_2+b_2+1
					\end{subarray}}
				\zeta(\alpha_0,\alpha_1,\cdots,\alpha_{a_1}+1,
					\beta_0,\beta_1,\cdots,\beta_{a_2}+1)+\\
				&+\sum_{
					\begin{subarray}{c}
						a_1+a_2+a_3=k\\
						b_1+b_2+b_3=\ell
					\end{subarray}}
				\mu_{\sigma(1)}^{a_1}
				\xi_{\sigma(1)}^{b_1}
				(\mu_{\sigma(1)}+\mu_{\sigma(2)})^{a_2}
				(\xi_{\sigma(1)}+\xi_{\sigma(2)})^{b_2}
				(
				\mu_{\sigma(1)}^{a_3}
				\xi_{\sigma(1)}^{b_3}+
				\mu_{\sigma(2)}^{a_3}
				\xi_{\sigma(2)}^{b_3}
				)\times\\
				&\quad\times\sum_{
					\begin{subarray}{c}
						|\boldsymbol\alpha|=a_1+b_1+1\\
						|\boldsymbol\beta|=a_2+b_2+1\\
						|\boldsymbol\gamma|=a_3+b_3+1
					\end{subarray}}
				\zeta(\alpha_0,\cdots,\alpha_{a_1},
					\beta_0,\beta_1,\cdots,\beta_{a_2}+
					\gamma_0,\gamma_1,\cdots,\gamma_{a_3}+1)\bigg],
			\end{align*}
			where $a_i$ and $b_i$ run over non-negative integers with the conditions that the sum of those are $k$ and $\ell$, and 
			$|\boldsymbol\alpha|:=\alpha_0+\cdots+\alpha_{a_1},
			|\boldsymbol\beta|:=\beta_0+\cdots+\beta_{a_2}$
			and
			$|\boldsymbol\gamma|:=\gamma_0+\cdots+\gamma_{a_3}$. 
			Therefore, for example $\displaystyle\sum_{
					\begin{subarray}{c}
						|\boldsymbol\alpha|=a_1+b_1+1\\
						|\boldsymbol\beta|=a_2+b_2+1
					\end{subarray}}$ means that $\alpha_0,\cdots,\alpha_{a_1}$ and
		$\beta_0,\cdots,\beta_{a_2}$ run over all 
		positive integers with $\alpha_0+\cdots+\alpha_{a_1}=a_1+b_1+1$,
		$\beta_0+\cdots+\beta_{a_2}=a_2+b_2+1$.
		\end{theorem}
		In {\bf Theorem \ref{theorem:sub theorem}}, when $k$ is an even number,
		putting $(\mu_1,\mu_2,\xi_1,\xi_2)=(1,-1,1,1)$ we can derive 
		the results of Eie, Liaw and Ong \cite[Main Theorem]{ELO} which is a
	 	generalization of the weighted sum formula of Ohno--Zudilin
	\cite[Theorem 3]{OZ}.
	We will derive \cite[Main Theorem]{ELO} from {\bf Theorem \ref{theorem:sub theorem}} in Section 3.

	The next one is the main theorem in this paper.

	\begin{theorem}
		\label{theorem:main theorem}
		For non-negative integers $k,\ell$ and six parameters (indeterminates) $\mu_1,\mu_2,\mu_3,\xi_1,\xi_2$ 
		and $\xi_3$, we have following relation:
	\begin{align*}
	&\sum_{
			\begin{subarray}{c}
						a_1+a_2+a_3=k\\
						b_1+b_2+b_3=\ell
			\end{subarray}}
		\mu_1^{a_1}\mu_2^{a_2}\mu_3^{a_3}\xi_1^{b_1}\xi_2^{b_2}\xi_3^{b_3}
		\zeta(a_1+b_1+2)\zeta(a_2+b_2+2)\zeta(a_3+b_3+2)\\
		=&\sum_{\sigma\in {\mathfrak S}_3}\bigg[
		\sum_{
			\begin{subarray}{c}
						a_1+a_2+a_3=k\\
						b_1+b_2+b_3=\ell
			\end{subarray}}
			P_{1,\sigma}
			\sum_{
			\begin{subarray}{c}
						|\boldsymbol\alpha|=a_1+b_1+1\\
						|\boldsymbol\beta|=a_2+b_2+1\\
						|\boldsymbol\gamma|=a_3+b_3+1
			\end{subarray}}
		\zeta(\alpha_0,\alpha_1,\cdots,\alpha_{a_1}+1,
			\beta_0,\beta_1,\cdots,\beta_{a_2}+1,
			\gamma_0,\gamma_1,\cdots,\gamma_{a_3}+1)+\\
			&+\sum_{
			\begin{subarray}{c}
						a_1+\cdots+a_4=k\\
						b_1+\cdots+b_4=\ell
			\end{subarray}}
			(P_{2,\sigma}+P_{3,\sigma})\times\\
			&\quad\times
			\sum_{
			\begin{subarray}{c}
						|\boldsymbol\alpha|=a_1+b_1+1\\
						|\boldsymbol\beta|=a_2+b_2+1\\
						|\boldsymbol\gamma|=a_3+b_3+1\\
						|\boldsymbol\delta|=a_4+b_4+1
			\end{subarray}}\hspace{-0.5cm}
		\zeta(\alpha_0,\cdots,\alpha_{a_1},
			\beta_0,\beta_1,\cdots,\beta_{a_2}+
			\gamma_0,\gamma_1,\cdots,\gamma_{a_3}+1,
			\delta_0,\delta_1,\cdots,\delta_{a_4}+1)+\\
			&+\sum_{
			\begin{subarray}{c}
						a_1+\cdots+a_5=k\\
						b_1+\cdots+b_5=\ell
			\end{subarray}}
			(P_{4,\sigma}+P_{5,\sigma}+P_{7,\sigma}+P_{12,\sigma})\times\\
			&\quad\times\sum_{
			\begin{subarray}{c}
						|\boldsymbol\alpha|=a_1+b_1+1\\
						|\boldsymbol\beta|=a_2+b_2+1\\
						|\boldsymbol\gamma|=a_3+b_3+1\\
						|\boldsymbol\delta|=a_4+b_4+1\\
						|\boldsymbol\varepsilon|=a_5+b_5+1
			\end{subarray}}\hspace{-0.5cm}
		\zeta(\alpha_0,\cdots,\alpha_{a_1},
			\beta_0,\beta_1,\cdots,\beta_{a_2}+
			\gamma_0,\gamma_1,\cdots,\gamma_{a_3},
			\delta_0,\delta_1,\cdots,\delta_{a_4}+\varepsilon_0,
			\varepsilon_1,\cdots,\varepsilon_{a_5}+1)+\\
			&+\sum_{
			\begin{subarray}{c}
						a_1+\cdots+a_4=k\\
						b_1+\cdots+b_4=\ell
			\end{subarray}}
			(P_{6,\sigma}+P_{11,\sigma})\times\\
			&\quad\times\sum_{
			\begin{subarray}{c}
						|\boldsymbol\alpha|=a_1+b_1+1\\
						|\boldsymbol\beta|=a_2+b_2+1\\
						|\boldsymbol\gamma|=a_3+b_3+1\\
						|\boldsymbol\delta|=a_4+b_4+1
			\end{subarray}}
		\zeta(\alpha_0,\alpha_1,\cdots,\alpha_{a_1}+1,
			\beta_0,\cdots,\beta_{a_2},
			\gamma_0,\gamma_1,\cdots,\gamma_{a_3}+
			\delta_0,\delta_1,\cdots,\delta_{a_4}+1)+\\
			&+\sum_{
			\begin{subarray}{c}
						a_1+\cdots+a_5=k\\
						b_1+\cdots+b_5=\ell
			\end{subarray}}
			(P_{8,\sigma}+P_{9,\sigma}+P_{10,\sigma}+P_{13,\sigma}+P_{14,\sigma}+P_{15,\sigma})\times\\
			&\quad\times\sum_{
			\begin{subarray}{c}
						|\boldsymbol\alpha|=a_1+b_1+1\\
						|\boldsymbol\beta|=a_2+b_2+1\\
						|\boldsymbol\gamma|=a_3+b_3+1\\
						|\boldsymbol\delta|=a_4+b_4+1\\
						|\boldsymbol\varepsilon|=a_5+b_5+1
			\end{subarray}}
		\zeta(\alpha_0,\cdots,\alpha_{a_1},
			\beta_0,\cdots,\beta_{a_2},
			\gamma_0,\gamma_1,\cdots,\gamma_{a_3}+
			\delta_0,\delta_1,\cdots,\delta_{a_4}+
			\varepsilon_0,\varepsilon_1,\cdots,\varepsilon_{a_5}+1)\bigg],
		\end{align*}
			where
		\begin{align*}
			P_{1,\sigma}=&
			\mu_{\sigma(1)}^{a_1}\mu_{\sigma(2)}^{a_2}\mu_{\sigma(3)}^{a_3}
			\xi_{\sigma(1)}^{b_1}\xi_{\sigma(2)}^{b_2}\xi_{\sigma(3)}^{b_3},\\
			P_{2,\sigma}
			=&\mu_{\sigma(1)}^{a_1}(\mu_{\sigma(1)}+\mu_{\sigma(2)})^{a_2}\mu_{\sigma(2)}^{a_3}\mu_{\sigma(3)}^{a_4}
			\xi_{\sigma(1)}^{b_1}(\xi_{\sigma(1)}+\xi_{\sigma(2)})^{b_2}\xi_{\sigma(2)}^{b_3}\xi_{\sigma(3)}^{b_4},\\
			P_{3,\sigma}
			=&\mu_{\sigma(1)}^{a_1+a_3}(\mu_{\sigma(1)}+\mu_{\sigma(2)})^{a_2}\mu_{\sigma(3)}^{a_4}
			\xi_{\sigma(1)}^{b_1+b_3}(\xi_{\sigma(1)}+\xi_{\sigma(2)})^{b_2}\xi_{\sigma(3)}^{b_4},\\
			P_{4,\sigma}
			=&\mu_{\sigma(1)}^{a_1+a_3}(\mu_{\sigma(1)}+\mu_{\sigma(2)})^{a_2}(\mu_{\sigma(1)}+\mu_{\sigma(3)})^{a_4}\mu_{\sigma(3)}^{a_5}
			\xi_{\sigma(1)}^{b_1+b_3}(\xi_{\sigma(1)}+\xi_{\sigma(2)})^{b_2}(\xi_{\sigma(1)}+\xi_{\sigma(3)})^{b_4}\xi_{\sigma(3)}^{b_5},\\
			P_{5,\sigma}
			=&\mu_{\sigma(1)}^{a_1+a_3+a_5}(\mu_{\sigma(1)}+\mu_{\sigma(2)})^{a_2}(\mu_{\sigma(1)}+\mu_{\sigma(3)})^{a_4}
			\xi_{\sigma(1)}^{b_1+b_3+b_5}(\xi_{\sigma(1)}+\xi_{\sigma(2)})^{b_2}(\xi_{\sigma(1)}+\xi_{\sigma(3)})^{b_4},\\
			P_{6,\sigma}
			=&\mu_{\sigma(1)}^{a_1}\mu_{\sigma(2)}^{a_2}(\mu_{\sigma(2)}+\mu_{\sigma(3)})^{a_3}\mu_{\sigma(3)}^{a_4}
			\xi_{\sigma(1)}^{b_1}\xi_{\sigma(2)}^{b_2}(\xi_{\sigma(2)}+\xi_{\sigma(3)})^{b_3}\xi_{\sigma(3)}^{b_4},\\
			P_{7,\sigma}
			=&\mu_{\sigma(1)}^{a_1}(\mu_{\sigma(1)}+\mu_{\sigma(2)})^{a_2}\mu_{\sigma(2)}^{a_3}(\mu_{\sigma(2)}+\mu_{\sigma(3)})^{a_4}\mu_{\sigma(3)}^{a_5}
			\xi_{\sigma(1)}^{b_1}(\xi_{\sigma(1)}+\xi_{\sigma(2)})^{b_2}\xi_{\sigma(2)}^{b_3}(\xi_{\sigma(2)}+\xi_{\sigma(3)})^{b_4}\xi_{\sigma(3)}^{b_5},\\
			P_{8,\sigma}
			=&\mu_{\sigma(1)}^{a_1}(\mu_{\sigma(1)}+\mu_{\sigma(2)})^{a_2}(\mu_{\sigma(1)}+\mu_{\sigma(2)}+\mu_{\sigma(3)})^{a_3}(\mu_{\sigma(2)}+\mu_{\sigma(3)})^{a_4}\mu_{\sigma(3)}^{a_5}\times\\
			&\times\xi_{\sigma(1)}^{b_1}(\xi_{\sigma(1)}+\xi_{\sigma(2)})^{b_2}(\xi_{\sigma(1)}+\xi_{\sigma(2)}+\xi_{\sigma(3)})^{b_3}(\xi_{\sigma(2)}+\xi_{\sigma(3)})^{b_4}\xi_{\sigma(3)}^{b_5},\\
			P_{9,\sigma}
			=&\mu_{\sigma(1)}^{a_1}(\mu_{\sigma(1)}+\mu_{\sigma(2)})^{a_2}(\mu_{\sigma(1)}+\mu_{\sigma(2)}+\mu_{\sigma(3)})^{a_3}(\mu_{\sigma(1)}+\mu_{\sigma(3)})^{a_4}\mu_{\sigma(3)}^{a_5}\times\\
			&\times\xi_{\sigma(1)}^{b_1}(\xi_{\sigma(1)}+\xi_{\sigma(2)})^{b_2}(\xi_{\sigma(1)}+\xi_{\sigma(2)}+\xi_{\sigma(3)})^{b_3}(\xi_{\sigma(1)}+\xi_{\sigma(3)})^{b_4}\xi_{\sigma(3)}^{b_5},\\
			P_{10,\sigma}
			=&\mu_{\sigma(1)}^{a_1+a_5}(\mu_{\sigma(1)}+\mu_{\sigma(2)})^{a_2}(\mu_{\sigma(1)}+\mu_{\sigma(2)}+\mu_{\sigma(3)})^{a_3}(\mu_{\sigma(1)}+\mu_{\sigma(3)})^{a_4}\times\\
			&\times\xi_{\sigma(1)}^{b_1+b_5}(\xi_{\sigma(1)}+\xi_{\sigma(2)})^{b_2}(\xi_{\sigma(1)}+\xi_{\sigma(2)}+\xi_{\sigma(3)})^{b_3}(\xi_{\sigma(1)}+\xi_{\sigma(3)})^{b_4},\\
			P_{11,\sigma}
			=&\mu_{\sigma(1)}^{a_1}\mu_{\sigma(2)}^{a_2+a_4}(\mu_{\sigma(2)}+\mu_{\sigma(3)})^{a_3}
			\xi_{\sigma(1)}^{b_1}\xi_{\sigma(2)}^{b_2+b_4}(\xi_{\sigma(2)}+\xi_{\sigma(3)})^{b_3},\\
			P_{12,\sigma}
			=&\mu_{\sigma(1)}^{a_1}(\mu_{\sigma(1)}+\mu_{\sigma(2)})^{a_2}\mu_{\sigma(2)}^{a_3+a_5}(\mu_{\sigma(2)}+\mu_{\sigma(3)})^{a_4}
			\xi_{\sigma(1)}^{b_1}(\xi_{\sigma(1)}+\xi_{\sigma(2)})^{b_2}\xi_{\sigma(2)}^{b_3+b_5}(\xi_{\sigma(2)}+\xi_{\sigma(3)})^{b_4},\\
			P_{13,\sigma}
			=&\mu_{\sigma(1)}^{a_1}(\mu_{\sigma(1)}+\mu_{\sigma(2)})^{a_2}(\mu_{\sigma(1)}+\mu_{\sigma(2)}+\mu_{\sigma(3)})^{a_3}(\mu_{\sigma(2)}+\mu_{\sigma(3)})^{a_4}\mu_{\sigma(2)}^{a_5}\times\\
			&\times\xi_{\sigma(1)}^{b_1}(\xi_{\sigma(1)}+\xi_{\sigma(2)})^{b_2}(\xi_{\sigma(1)}+\xi_{\sigma(2)}+\xi_{\sigma(3)})^{b_3}(\xi_{\sigma(2)}+\xi_{\sigma(3)})^{b_4}\xi_{\sigma(2)}^{b_5},\\
			P_{14,\sigma}
			=&\mu_{\sigma(1)}^{a_1}(\mu_{\sigma(1)}+\mu_{\sigma(2)})^{a_2+a_4}(\mu_{\sigma(1)}+\mu_{\sigma(2)}+\mu_{\sigma(3)})^{a_3}\mu_{\sigma(2)}^{a_5}\times\\
			&\times\xi_{\sigma(1)}^{b_1}(\xi_{\sigma(1)}+\xi_{\sigma(2)})^{b_2+b_4}(\xi_{\sigma(1)}+\xi_{\sigma(2)}+\xi_{\sigma(3)})^{b_3}\xi_{\sigma(2)}^{b_5},\\
			P_{15,\sigma}
			=&\mu_{\sigma(1)}^{a_1+a_5}(\mu_{\sigma(1)}+\mu_{\sigma(2)})^{a_2+a_4}(\mu_{\sigma(1)}+\mu_{\sigma(2)}+\mu_{\sigma(3)})^{a_3}
			\xi_{\sigma(1)}^{b_1+b_5}(\xi_{\sigma(1)}+\xi_{\sigma(2)})^{b_2+b_4}(\xi_{\sigma(1)}+\xi_{\sigma(2)}+\xi_{\sigma(3)})^{b_3}.
		\end{align*}
	\end{theorem}

	\begin{remark}
		The formulas of  {\bf Theorem \ref{theorem:sub theorem}} and
		 {\bf Theorem \ref{theorem:main theorem}} have several parameters,
		so one may get new relations by comparing the parameters of the both sides
		or differentiating partially with respect to parameters.
	\end{remark}

\section{Proof of Theorem \ref{theorem:sub theorem} and \ref{theorem:main theorem}}
	In this section, we prove {\bf Theorem \ref{theorem:main theorem}}.
	The basic structure of the proof is the same as in \cite{ELO}.
	For $k, \ell\in\mathbb{Z}_{\geq0}$, and 6 parameters $\mu_1, \mu_2, \mu_3, \xi_1, \xi_2$ and $\xi_3$,
	we define $I_{k,\ell}(\mu_1,\mu_2,\mu_3,\xi_1,\xi_2,\xi_3)$ as the following integral:
		\begin{align*}
	 		&I_{k,\ell}(\mu_1,\mu_2,\mu_3,\xi_1,\xi_2,\xi_3)\\
			&:=\frac{1}{k!\ell!}\int\limits_{
					\begin{subarray}{c}
						0<s_1<s_2<1\\
						0<t_1<t_2<1\\
						0<u_1<u_2<1
					\end{subarray}}
			\left(\mu_1\log\frac{1-s_1}{1-s_2}
			+\mu_2\log\frac{1-t_1}{1-t_2}
			+\mu_3\log\frac{1-u_1}{1-u_2}\right)^k\times\\
			&\times\left(\xi_1\log\frac{s_2}{s_1}
			+\xi_2\log\frac{t_2}{t_1}
			+\xi_3\log\frac{u_2}{u_1}\right)^{\ell}
			\frac{ds_1ds_2dt_1dt_2du_1du_2}{(1-s_1)s_2(1-t_1)t_2(1-u_1)u_2}.
	\end{align*}
	To obtain {\bf Theorem \ref{theorem:main theorem}}, we calculate this integral in two ways.\\
	$\bullet$ The first calculations\\
	First, we expand the factors of the integrand simply as follows:
		\begin{align*}
			&\left(\mu_1\log\frac{1-s_1}{1-s_2}
			+\mu_2\log\frac{1-t_1}{1-t_2}
			+\mu_3\log\frac{1-u_1}{1-u_2}\right)^k\\
			=&\sum_{a_1+a_2+a_3=k}\frac{k!}{a_1!a_2!a_3!}
			\mu_1^{a_1}\mu_2^{a_2}\mu_3^{a_3}
			\left(\log\frac{1-s_1}{1-s_2}\right)^{a_1}
			\left(\log\frac{1-t_1}{1-t_2}\right)^{a_2}
			\left(\log\frac{1-u_1}{1-u_2}\right)^{a_3},\\
			&\left(\xi_1\log\frac{s_2}{s_1}
			+\xi_2\log\frac{t_2}{t_1}
			+\xi_3\log\frac{u_2}{u_1}\right)^{\ell}\\
			=&\sum_{b_1+b_2+b_3=\ell}\frac{\ell!}{b_1!b_2!b_3!}
			\xi_1^{b_1}\xi_2^{b_2}\xi_3^{b_3}
			\left(\log\frac{s_2}{s_1}\right)^{b_1}
			\left(\log\frac{t_2}{t_1}\right)^{b_2}
			\left(\log\frac{u_2}{u_1}\right)^{b_3}.
		\end{align*}
		Substituting these expansions, we get
		\begin{align*}
			&I_{k,\ell}(\mu_1,\mu_2,\mu_3,\xi_1,\xi_2,\xi_3)\\
				=&\sum_{
					\begin{subarray}{c}
						a_1+a_2+a_3=k\\
						b_1+b_2+b_3=\ell
					\end{subarray}}
			\mu_1^{a_1}\mu_2^{a_2}\mu_3^{a_3}\xi_1^{b_1}\xi_2^{b_2}\xi_3^{b_3}\times\\
			&\times\left\{\frac{1}{a_1!b_1!}\int\limits_{0<s_1<s_2<1}
			\left(\log\frac{1-s_1}{1-s_2}\right)^{a_1}
			\left(\log\frac{s_2}{s_1}\right)^{b_1}
			\frac{ds_1ds_2}{(1-s_1)s_2}
			\right\}\times\\
			&\times
			\left\{\frac{1}{a_2!b_2!}\int\limits_{0<t_1<t_2<1}
			\left(\log\frac{1-t_1}{1-t_2}\right)^{a_2}
			\left(\log\frac{t_2}{t_1}\right)^{b_2}
			\frac{dt_1dt_2}{(1-t_1)t_2}
			\right\}\times\\
			&\times\left\{\frac{1}{a_3!b_3!}\int\limits_{0<u_1<u_2<1}
			\left(\log\frac{1-u_1}{1-u_2}\right)^{a_3}
			\left(\log\frac{u_2}{u_1}\right)^{b_3}
			\frac{du_1du_2}{(1-u_1)u_2}
			\right\}.
		\end{align*}

	Here let us consider $\displaystyle\left(\log\frac{1-s_1}{1-s_2}\right)^a$.
	We can rewrite it as
		\begin{align*}
			\displaystyle
					\left(\log\frac{1-s_1}{1-s_2}\right)^a
					=&\left(\int_{s_1}^{s_2}\frac{dp}{1-p}\right)^a\\
					=&\int\limits_{
						\begin{subarray}{c}
							s_1<p_1<s_2\\
							\vdots\\
							s_1<p_a<s_2
						\end{subarray}}
				\prod_{i=1}^{a}\frac{dp_i}{1-p_i}\\
				=&\sum_{\sigma\in {\mathfrak S}_a}\int\limits_{s_1<p_{\sigma(1)}<\cdots<p_{\sigma(a)}<s_2}
				\prod_{i=1}^{a}\frac{dp_i}{1-p_i}\\
				=&a!\int\limits_{s_1<p_{1}<\cdots<p_{a}<s_2}
				\prod_{i=1}^{a}\frac{dp_i}{1-p_i}.
		\end{align*} 

	We can rewrite $\displaystyle\left(\log\frac{s_2}{s_1}\right)^b$ similarly, so we get the  following lemma.
		\begin{lemma}
			\label{lemma:repeat integral}
			For $a,b\in\mathbb{Z}_{\geq0}$, and $0<s_1<s_2<1$ we have
				\[
					\displaystyle
					\left(\log\frac{1-s_1}{1-s_2}\right)^a
					=a!\int\limits_{s_1<p_1<\cdots<p_a<s_2}
					\prod_{i=1}^{a}\frac{dp_i}{1-p_i}
				\]
			and
				\[
					\displaystyle\left(\log\frac{s_2}{s_1}\right)^b
					=b!\int\limits_{s_1<q_1<\cdots<q_b<s_2}
					\prod_{j=1}^{b}\frac{dq_j}{q_j},
				\]
		where the empty product is to be understood as {\rm 1}.
		\end{lemma}
	Using the above lemma, we get
		\begin{equation}
			\label{eq:1}
			\begin{split}
			&\frac{1}{a_1!b_1!}\int\limits_{0<s_1<s_2<1}
			\left(\log\frac{1-s_1}{1-s_2}\right)^{a_1}
			\left(\log\frac{s_2}{s_1}\right)^{b_1}
			\frac{ds_1ds_2}{(1-s_1)s_2}\\
			=&\int\limits_{
					\begin{subarray}{c}
						0<s_1<s_2<1\\
						s_1<p_{1}<\cdots<p_{a_1}<s_2\\
						s_1<q_{1}<\cdots<q_{b_1}<s_2
					\end{subarray}}
				\frac{ds_1}{1-s_1}
				\prod_{i=1}^{a_1}\frac{dp_{i}}{1-p_{i}}
				\prod_{j=1}^{b_1}\frac{dq_{j}}{q_{j}}
				\frac{ds_2}{s_2}\\
				=&\sum_{(r_1,\cdots,r_{a_1+b_1})}
				\int\limits_
				{\begin{subarray}{c}
				0<s_1<s_2<1\\
				s_1<r_1<\cdots<r_{a_1+b_1}<s_2
				\end{subarray}}
				\frac{ds_1}{1-s_1}
				\prod_{i=1}^{a_1}\frac{dp_{i}}{1-p_{i}}
				\prod_{j=1}^{b_1}\frac{dq_{j}}{q_{j}}
				\frac{ds_2}{s_2},
		\end{split}
		\end{equation}
	where the summation runs over all tuples $(r_1,\cdots,r_{a_1+b_1})$ such that 
	\[
		\{r_1,\cdots,r_{a_1+b_1}\}=\{p_1,\cdots,p_{a_1}\}\cup\{q_1,\cdots,q_{b_1}\}
	\]
	and 
	\[
		p_1<\cdots<p_{a_1},
		q_1<\cdots<q_{b_1}.
	\]
	Then each integral gives a multiple zeta value, which implies that the above is
		\begin{align*}
			=\sum_{
					\begin{subarray}{c}
						|\boldsymbol\alpha|=a_1+b_1+1
					\end{subarray}}
				\zeta(\alpha_0,\alpha_1,\cdots,\alpha_{a_1}+1).
		\end{align*}
	By using the sum formula \cite[Proposition]{Gr}, we have
		\begin{align*}
			I_{k,\ell}(\mu_1,\mu_2,\mu_3,\xi_1,\xi_2,\xi_3)
			=\sum_{
					\begin{subarray}{c}
						a_1+a_2+a_3=k\\
						b_1+b_2+b_3=\ell
					\end{subarray}}
			\mu_1^{a_1}\mu_2^{a_2}\mu_3^{a_3}\xi_1^{b_1}\xi_2^{b_2}\xi_3^{b_3}
			\zeta(a_1+b_1+2)\zeta(a_2+b_2+2)\zeta(a_3+b_3+2).
		\end{align*}
	This is the end of the first calculations.\\

	$\bullet$ The second calculations\\
	We divide the region $0<s_1<s_2<1, 0<t_1<t_2<1, 0<u_1<u_2<1$ to ninety regions,
	according to the order of magnitude of variables, and calculate the integral on each regions.
	But, as we will see just below, because of the symmetry of variables $s,t$ and $u$, it is sufficient to calculate on the following fifteen regions:
		\begin{align*}
			&D_1:0<s_1<s_2<t_1<t_2<u_1<u_2<1,\quad D_2:0<s_1<t_1<s_2<t_2<u_1<u_2<1,\\
			&D_3:0<s_1<t_1<t_2<s_2<u_1<u_2<1,\quad D_4:0<s_1<t_1<t_2<u_1<s_2<u_2<1,\\
			&D_5:0<s_1<t_1<t_2<u_1<u_2<s_2<1,\\
			&D_6:0<s_1<s_2<t_1<u_1<t_2<u_2<1,\quad D_7:0<s_1<t_1<s_2<u_1<t_2<u_2<1,\\
			&D_8:0<s_1<t_1<u_1<s_2<t_2<u_2<1,\quad D_9:0<s_1<t_1<u_1<t_2<s_2<u_2<1,\\
			&D_{10}:0<s_1<t_1<u_1<t_2<u_2<s_2<1,\\
			&D_{11}:0<s_1<s_2<t_1<u_1<u_2<t_2<1,\quad D_{12}:0<s_1<t_1<s_2<u_1<u_2<t_2<1,\\
			&D_{13}:0<s_1<t_1<u_1<s_2<u_2<t_2<1,\quad D_{14}:0<s_1<t_1<u_1<u_2<s_2<t_2<1,\\
			&D_{15}:0<s_1<t_1<u_1<u_2<t_2<s_2<1.
		\end{align*}
	We define $I_{k,\ell,m}(\mu_1,\mu_2,\mu_3,\xi_1,\xi_2,\xi_3)$ as the integral on the region $D_{m}$ ($m=1,\cdots,15$):
		\begin{align*}
	 		&I_{k,\ell,m}(\mu_1,\mu_2,\mu_3,\xi_1,\xi_2,\xi_3)\\
			&:=\frac{1}{k!\ell!}\int\limits_{D_{m}}
			\left(\mu_1\log\frac{1-s_1}{1-s_2}
			+\mu_2\log\frac{1-t_1}{1-t_2}
			+\mu_3\log\frac{1-u_1}{1-u_2}\right)^k\times\\
			&\times\left(\xi_1\log\frac{s_2}{s_1}
			+\xi_2\log\frac{t_2}{t_1}
			+\xi_3\log\frac{u_2}{u_1}\right)^{\ell}
			\frac{ds_1ds_2dt_1dt_2du_1du_2}{(1-s_1)s_2(1-t_1)t_2(1-u_1)u_2}.
	\end{align*}
	How to treat the remaining seventy five regions?
	For example, we can see that the integration on the region $0<u_1<u_2<s_1<s_2<t_1<t_2<1$ 
	which is one of the remaining seventy five regions is written by using 
	the integration on $D_{1}$ as follows:

	\begin{align*}
		&\frac{1}{k!\ell!}\int\limits_{0<u_1<u_2<s_1<s_2<t_1<t_2<1}
		\left(\mu_1\log\frac{1-s_1}{1-s_2}
		+\mu_2\log\frac{1-t_1}{1-t_2}
		+\mu_3\log\frac{1-u_1}{1-u_2}\right)^k\times\\
		&\times\left(\xi_1\log\frac{s_2}{s_1}
		+\xi_2\log\frac{t_2}{t_1}
		+\xi_3\log\frac{u_2}{u_1}\right)^{\ell}
		\frac{ds_1ds_2dt_1dt_2du_1du_2}{(1-s_1)s_2(1-t_1)t_2(1-u_1)u_2}\\
		=&I_{k,\ell,1}(\mu_3,\mu_1,\mu_2,\xi_3,\xi_1,\xi_2).
	\end{align*}

	Therefore the integrals on the remaining seventy five regions are obtained by changing parameters
	in the integrals on $D_{1}, D_{2},\cdots,D_{14}$ or $D_{15}$.
	After all, we obtain
		\begin{align*}
			I_{k,\ell}(\mu_1,\mu_2,\mu_3,\xi_1,\xi_2,\xi_3)
			=\sum_{\sigma\in {\mathfrak S}_{3}}\sum_{m=1}^{15} I_{k,\ell,m}(\mu_{\sigma(1)},\mu_{\sigma(2)},\mu_{\sigma(3)},\xi_{\sigma(1)},\xi_{\sigma(2)},\xi_{\sigma(3)}).
		\end{align*}  

	To save pages, we present only the calculations on the region $D_{4}:0<s_1<t_1<t_2<u_1<s_2<u_2<1$ in this paper.
	We want to write the integral as an explicit sum of MZVs, hence we need to
	modify some terms in the integrand before we expand the integrand.
	Those modifications are the most important point in this proof.
	We expand the integrand as follows:
		\begin{equation}
				\label{eq:2}
				\begin{split}
			&\left(\mu_{\sigma(1)}\log\frac{1-s_1}{1-s_2}
			+\mu_{\sigma(2)}\log\frac{1-t_1}{1-t_2}
			+\mu_{\sigma(3)}\log\frac{1-u_1}{1-u_2}\right)^k\\
				=&\left\{
			\mu_{\sigma(1)}\log\frac{1-s_1}{1-t_1}
			+(\mu_{\sigma(1)}+\mu_{\sigma(2)})\log\frac{1-t_1}{1-t_2}
			+\mu_{\sigma(1)}\log\frac{1-t_2}{1-u_1}+\right.\\
			&\quad\left.+(\mu_{\sigma(1)}+\mu_{\sigma(3)})\log\frac{1-u_1}{1-s_2}
			+\mu_{\sigma(3)}\log\frac{1-s_2}{1-u_2}
			\right\}^k\\
			=&\sum_{a_1+a_2+a_3+a_4+a_5=k}\frac{k!}{a_1!a_2!a_3!a_4!a_5!}
			\mu_{\sigma(1)}^{a_1+a_3}(\mu_{\sigma(1)}+\mu_{\sigma(2)})^{a_2}(\mu_{\sigma(1)}+\mu_{\sigma(3)})^{a_4}\mu_{\sigma(3)}^{a_5}\times\\
			&\times\left(\log\frac{1-s_1}{1-t_1}\right)^{a_1}
			\left(\log\frac{1-t_1}{1-t_2}\right)^{a_2}
			\left(\log\frac{1-t_2}{1-u_1}\right)^{a_3}
			\left(\log\frac{1-u_1}{1-s_2}\right)^{a_4}
			\left(\log\frac{1-s_2}{1-u_2}\right)^{a_5},
			\end{split}
			\end{equation}
			and
			\begin{equation}
				\label{eq:3}
				\begin{split}
			&\left(\xi_{\sigma(1)}\log\frac{s_2}{s_1}
			+\xi_{\sigma(2)}\log\frac{t_2}{t_1}
			+\xi_{\sigma(3)}\log\frac{u_2}{u_1}\right)^{\ell}\\
			=&\left\{
			\xi_{\sigma(1)}\log\frac{t_1}{s_1}
			+(\xi_{\sigma(1)}+\xi_{\sigma(2)})\log\frac{t_2}{t_1}
			+\xi_{\sigma(1)}\log\frac{u_1}{t_2}+\right.\\
			&\quad\left.+(\xi_{\sigma(1)}+\xi_{\sigma(3)})\log\frac{s_2}{u_1}
			+\xi_{\sigma(3)}\log\frac{u_2}{s_2}
			\right\}^{\ell}\\
			=&\sum_{b_1+b_2+b_3+b_4+b_5=\ell}\frac{\ell!}{b_1!b_2!b_3!b_4!b_5!}
			\xi_{\sigma(1)}^{b_1+b_3}(\xi_{\sigma(1)}+\xi_{\sigma(2)})^{b_2}(\xi_{\sigma(1)}+\xi_{\sigma(3)})^{b_4}\xi_{\sigma(3)}^{b_5}\times\\
			&\times\left(\log\frac{t_1}{s_1}\right)^{b_1}
			\left(\log\frac{t_2}{t_1}\right)^{b_2}
			\left(\log\frac{u_1}{t_2}\right)^{b_3}
			\left(\log\frac{s_2}{u_1}\right)^{b_4}
			\left(\log\frac{u_2}{s_2}\right)^{b_5}.
		\end{split}
		\end{equation}

	The above modifications are based on the following observation.
	First consider $\log\frac{1-s_1}{1-s_2}$.
	In this case, there are $t_1,t_2$ and $u_1$ between $s_1$ and $s_2$, so we modify 
	$\log\frac{1-s_1}{1-s_2}$ to the following form.
		\begin{align*}
			\log\frac{1-s_1}{1-s_2}
			&=\int_{s_1}^{s_2}\frac{dp}{1-p}\\
			&=\int_{s_1}^{t_1}+\int_{t_1}^{t_2}+\int_{t_2}^{u_1}+\int_{u_1}^{s_2}\frac{dp}{1-p}\\
			&=\log\frac{1-s_1}{1-t_1}+\log\frac{1-t_1}{1-t_2}+\log\frac{1-t_2}{1-u_1}+\log\frac{1-u_1}{1-s_2}.
		\end{align*}
	There is $s_2$ between $u_1$ and $u_2$, so we modify $\log\frac{1-u_1}{1-u_2}$ similarly:
		\begin{align*}
		\log\frac{1-u_1}{1-u_2}
			&=\int_{u_1}^{u_2}\frac{dp}{1-p}\\
			&=\int_{u_1}^{s_2}+\int_{s_2}^{u_2}\frac{dp}{1-p}\\
			&=\log\frac{1-u_1}{1-s_2}+\log\frac{1-s_2}{1-u_2}.
		\end{align*}
	These modifications give (\ref{eq:2}), and similarly we can show (\ref{eq:3}).
	Substituting the above modified expansions and using 
	{\bf Lemma \ref{lemma:repeat integral}} and arguing in the same way as (\ref{eq:1}),
	we get
		\begin{align*}
			 &I_{k,\ell,4}(\mu_{\sigma(1)},\mu_{\sigma(2)},\mu_{\sigma(3)},\xi_{\sigma(1)},\xi_{\sigma(2)},\xi_{\sigma(3)})\\
			=&\sum_{
			\begin{subarray}{c}
						a_1+\cdots+a_5=k\\
						b_1+\cdots+b_5=\ell
			\end{subarray}}
		P_{4,\sigma}
		\int\limits_{D_4^{\prime}}
		\frac{ds_1}{1-s_1}
		\prod_{i=1}^{a_1}\frac{dp_{i}}{1-p_{i}}
		\prod_{j=1}^{b_1}\frac{dq_{j}}{q_{j}}
		\frac{dt_1}{1-t_1}
		\prod_{i=a_1+1}^{a_1+a_2}\frac{dp_{i}}{1-p_{i}}
		\prod_{j=b_1+1}^{b_1+b_2}\frac{dq_{j}}{q_{j}}
		\frac{dt_2}{t_2}\times\\
		&\times\prod_{i=a_1+a_2+1}^{a_1+a_2+a_3}\frac{dp_{i}}{1-p_{i}}
		\prod_{j=b_1+b_2+1}^{b_1+b_2+b_3}\frac{dq_{j}}{q_{j}}
		\frac{du_1}{1-u_1}
		\prod_{i=a_1+a_2+a_3+1}^{a_1+\cdots+a_4}\frac{dp_{i}}{1-p_{i}}
		\prod_{j=b_1+b_2+b_3+1}^{b_1+\cdots+b_4}\frac{dq_{j}}{q_{j}}
		\frac{ds_2}{s_2}\times\\
		&\times\prod_{i=a_1+\cdots+a_4+1}^{a_1+\cdots+a_5}\frac{dp_{i}}{1-p_{i}}
		\prod_{j=b_1+\cdots+b_4+1}^{b_1+\cdots+b_5}\frac{dq_{j}}{q_{j}}
		\frac{du_2}{u_2}\\
		=&\sum_{
			\begin{subarray}{c}
						a_1+\cdots+a_5=k\\
						b_1+\cdots+b_5=\ell
			\end{subarray}}
			P_{4,\sigma}\sum_{
			\begin{subarray}{c}
						|\boldsymbol\alpha|=a_1+b_1+1\\
						|\boldsymbol\beta|=a_2+b_2+1\\
						|\boldsymbol\gamma|=a_3+b_3+1\\
						|\boldsymbol\delta|=a_4+b_4+1\\
						|\boldsymbol\varepsilon|=a_5+b_5+1
			\end{subarray}}
		\zeta(\alpha_0,\cdots,\alpha_{a_1},
			\beta_0,\beta_1,\cdots,\beta_{a_2}+
			\gamma_0,\gamma_1,\cdots,\gamma_{a_3},\\
			&\hspace{6.0cm}\quad\delta_0,\delta_1,\cdots,\delta_{a_4}+\varepsilon_0,
			\varepsilon_1,\cdots,\varepsilon_{a_5}+1),
		\end{align*}
	where
		\begin{align*}
			D_{4}^{\prime}
			=\left\{(s_1,s_2,t_1,t_2,u_1,u_2,p_1,\cdots,p_{k},q_1,\cdots,q_{\ell})\in (0,1)^{k+\ell+4}
			\quad\Bigg|
				\begin{subarray}{c}
						0<s_1<t_1<t_2<u_1<s_2<u_2<1\\
						s_1<p_{1}<\cdots<p_{a_1}<t_1\\
      					s_1<q_{1}<\cdots<q_{b_1}<t_1\\
						t_1<p_{a_1+1}<\cdots<p_{a_1+a_2}<t_2\\
      					t_1<q_{b_1+1}<\cdots<q_{b_1+b_2}<t_2\\
      					t_2<p_{a_1+a_2+1}<\cdots<p_{a_1+a_2+a_3}<u_1\\
      					t_2<q_{b_1+b_2+1}<\cdots<q_{b_1+b_2+b_3}<u_1\\    			
      					u_1<p_{a_1+a_2+a_3+1}<\cdots<p_{a_1+\cdots+a_4}<s_2\\
      					u_1<q_{b_1+b_2+b_3+1}<\cdots<q_{b_1+\cdots+b_4}<s_2\\
      					s_2<p_{a_1+\cdots+a_4+1}<\cdots<p_{a_1+\cdots+a_5}<u_2\\
      					s_2<q_{b_1+\cdots+b_4+1}<\cdots<q_{b_1+\cdots+b_5}<u_2
				\end{subarray}
			\right\}
		\end{align*}
	and 
		\begin{align*}
			P_{4,\sigma}
			=\mu_{\sigma(1)}^{a_1+a_3}(\mu_{\sigma(1)}+\mu_{\sigma(2)})^{a_2}(\mu_{\sigma(1)}+\mu_{\sigma(3)})^{a_4}\mu_{\sigma(3)}^{a_5}
			\xi_{\sigma(1)}^{b_1+b_3}(\xi_{\sigma(1)}+\xi_{\sigma(2)})^{b_2}(\xi_{\sigma(1)}+\xi_{\sigma(3)})^{b_4}\xi_{\sigma(3)}^{b_5}.
		\end{align*}
	This is the end of the calculations for $I_{k,\ell,4}(\mu_{\sigma(1)},\mu_{\sigma(2)},\mu_{\sigma(3)},\xi_{\sigma(1)},\xi_{\sigma(2)},\xi_{\sigma(3)})$.
	When one calculate the other $I_{k,\ell,m}(\mu_{\sigma(1)},\mu_{\sigma(2)},\mu_{\sigma(3)},\xi_{\sigma(1)},\xi_{\sigma(2)},\xi_{\sigma(3)})$ by the same way of modifications,
	then the $P_{m,\sigma}$ which is in {\bf Theorem \ref{theorem:main theorem}} appears as the coefficient of modified integrand for the integral of $D_{m}$.
	Namely, one obtain 
	\[
		\displaystyle I_{k,\ell,m}(\mu_{\sigma(1)},\mu_{\sigma(2)},\mu_{\sigma(3)},\xi_{\sigma(1)},\xi_{\sigma(2)},\xi_{\sigma(3)})=\sum_{a,b} P_{m,\sigma}\times{(\rm the~sum~of~MZVs)}
	\] as above.
	Moreover, several types of the sum of MZVs will appear in the calculations, 
	but one can notice that these types of the sum of MZVs depend only on  
	the subscript of divided regions.
	Namely, the sum of MZVs which appears after the calculations 
	for $D_{2}$ and $D_{3}$ are the same, for $D_{4},D_{5},D_{7}$ and $D_{12}$ are the same,
	for $D_{6}$ and $D_{11}$ are the same, and for $D_{8},D_{9},D_{10},D_{13},D_{14}$ and $D_{15}$
	are the same.
	Noting this point, we obtain

		\begin{align*}
		&I_{k,\ell}(\mu_1,\mu_2,\mu_3,\xi_1,\xi_2,\xi_3)\\
		=&\sum_{\sigma\in {\mathfrak S}_3}\bigg[
		\sum_{
			\begin{subarray}{c}
						a_1+a_2+a_3=k\\
						b_1+b_2+b_3=\ell
			\end{subarray}}
			P_{1,\sigma}
			\sum_{
			\begin{subarray}{c}
						|\boldsymbol\alpha|=a_1+b_1+1\\
						|\boldsymbol\beta|=a_2+b_2+1\\
						|\boldsymbol\gamma|=a_3+b_3+1
			\end{subarray}}\hspace{-0.7cm}
		\zeta(\alpha_0,\alpha_1,\cdots,\alpha_{a_1}+1,
			\beta_0,\beta_1,\cdots,\beta_{a_2}+1,
			\gamma_0,\gamma_1,\cdots,\gamma_{a_3}+1)+\\
			&+\sum_{
			\begin{subarray}{c}
						a_1+\cdots+a_4=k\\
						b_1+\cdots+b_4=\ell
			\end{subarray}}
			(P_{2,\sigma}+P_{3,\sigma})\times\\
			&\quad\times
			\sum_{
			\begin{subarray}{c}
						|\boldsymbol\alpha|=a_1+b_1+1\\
						|\boldsymbol\beta|=a_2+b_2+1\\
						|\boldsymbol\gamma|=a_3+b_3+1\\
						|\boldsymbol\delta|=a_4+b_4+1
			\end{subarray}}
		\zeta(\alpha_0,\cdots,\alpha_{a_1},
			\beta_0,\beta_1,\cdots,\beta_{a_2}+
			\gamma_0,\gamma_1,\cdots,\gamma_{a_3}+1,
			\delta_0,\delta_1,\cdots,\delta_{a_4}+1)+\\
			&+\sum_{
			\begin{subarray}{c}
						a_1+\cdots+a_5=k\\
						b_1+\cdots+b_5=\ell
			\end{subarray}}
			(P_{4,\sigma}+P_{5,\sigma}+P_{7,\sigma}+P_{12,\sigma})\times\\
			&\quad\times\sum_{
			\begin{subarray}{c}
						|\boldsymbol\alpha|=a_1+b_1+1\\
						|\boldsymbol\beta|=a_2+b_2+1\\
						|\boldsymbol\gamma|=a_3+b_3+1\\
						|\boldsymbol\delta|=a_4+b_4+1\\
						|\boldsymbol\varepsilon|=a_5+b_5+1
			\end{subarray}}
		\zeta(\alpha_0,\cdots,\alpha_{a_1},
			\beta_0,\beta_1,\cdots,\beta_{a_2}+
			\gamma_0,\gamma_1,\cdots,\gamma_{a_3},
			\delta_0,\delta_1,\cdots,\delta_{a_4}+\varepsilon_0,
			\varepsilon_1,\cdots,\varepsilon_{a_5}+1)+\\
			&+\sum_{
			\begin{subarray}{c}
						a_1+\cdots+a_4=k\\
						b_1+\cdots+b_4=\ell
			\end{subarray}}
			(P_{6,\sigma}+P_{11,\sigma})\times\\
			&\quad\times\sum_{
			\begin{subarray}{c}
						|\boldsymbol\alpha|=a_1+b_1+1\\
						|\boldsymbol\beta|=a_2+b_2+1\\
						|\boldsymbol\gamma|=a_3+b_3+1\\
						|\boldsymbol\delta|=a_4+b_4+1
			\end{subarray}}
		\zeta(\alpha_0,\alpha_1,\cdots,\alpha_{a_1}+1,
			\beta_0,\cdots,\beta_{a_2},
			\gamma_0,\gamma_1,\cdots,\gamma_{a_3}+
			\delta_0,\delta_1,\cdots,\delta_{a_4}+1)+\\
			&+\sum_{
			\begin{subarray}{c}
						a_1+\cdots+a_5=k\\
						b_1+\cdots+b_5=\ell
			\end{subarray}}
			(P_{8,\sigma}+P_{9,\sigma}+P_{10,\sigma}+P_{13,\sigma}+P_{14,\sigma}+P_{15,\sigma})\times\\
			&\quad\times\sum_{
			\begin{subarray}{c}
						|\boldsymbol\alpha|=a_1+b_1+1\\
						|\boldsymbol\beta|=a_2+b_2+1\\
						|\boldsymbol\gamma|=a_3+b_3+1\\
						|\boldsymbol\delta|=a_4+b_4+1\\
						|\boldsymbol\varepsilon|=a_5+b_5+1
			\end{subarray}}
		\zeta(\alpha_0,\cdots,\alpha_{a_1},
			\beta_0,\cdots,\beta_{a_2},
			\gamma_0,\gamma_1,\cdots,\gamma_{a_3}+
			\delta_0,\delta_1,\cdots,\delta_{a_4}+
			\varepsilon_0,\varepsilon_1,\cdots,\varepsilon_{a_5}+1)\bigg].
		\end{align*}
	Combining the conclusions of the first calculations and the second calculations,
	we obtain the asserted relations of  {\bf Theorem \ref{theorem:main theorem}}. 

		To prove {\bf Theorem \ref{theorem:sub theorem}}, it is sufficient to calculate similarly the following integrals:
		\begin{align*}
	 		&I_{k,\ell}(\mu_1,\mu_2,\xi_1,\xi_2)\\
			&:=\frac{1}{k!\ell!}\int\limits_{
					\begin{subarray}{c}
						0<s_1<s_2<1\\
						0<t_1<t_2<1
					\end{subarray}}
			\left(\mu_1\log\frac{1-s_1}{1-s_2}
			+\mu_2\log\frac{1-t_1}{1-t_2}\right)^k\times\\
			&\times\left(\xi_1\log\frac{s_2}{s_1}
			+\xi_2\log\frac{t_2}{t_1}\right)^{\ell}
			\frac{ds_1ds_2dt_1dt_2}{(1-s_1)s_2(1-t_1)t_2}.
	\end{align*}
	We omit the details of the proof.

\section{Deduction of the result of Eie, Liaw and Ong, and some more generalizations}
	Let us consider a special case of {\bf Theorem \ref{theorem:sub theorem}}.
	Putting $\xi_{1}=\xi_{2}=\xi\neq0$ in the formula of {\bf Theorem \ref{theorem:sub theorem}},
	and using the harmonic product $\zeta(a)\zeta(b)=\zeta(a,b)+\zeta(b,a)+\zeta(a+b)$
	for the left hand side, we find that $\xi^{\ell}$ parts of the both sides are cancelled 
	with each other, and
		\begin{equation}
			\label{eq:after}
			\begin{split}
			&(\ell+1)\zeta(k+\ell+4)
			\sum_{a_1+a_2=k}
			\mu_{1}^{a_1}\mu_{2}^{a_2}	
			+\sum_{
					\begin{subarray}{c}
						a_1+a_2=k\\
						b_1+b_2=\ell
					\end{subarray}}
			(\mu_{1}^{a_1}\mu_{2}^{a_2}+\mu_{2}^{a_1}\mu_{1}^{a_2})
			\zeta(a_2+b_2+2,a_1+b_1+2)\\
				=&
				\sum_{
					\begin{subarray}{c}
						a_1+a_2=k\\
						b_1+b_2=\ell
					\end{subarray}}
				(\mu_{1}^{a_1}\mu_{2}^{a_2}+\mu_{2}^{a_1}\mu_{1}^{a_2})
				\sum_{
					\begin{subarray}{c}
						|\boldsymbol\alpha|=a_1+b_1+1\\
						|\boldsymbol\beta|=a_2+b_2+1
					\end{subarray}}
				\zeta(\alpha_0,\alpha_1,\cdots,\alpha_{a_1}+1,
					\beta_0,\beta_1,\cdots,\beta_{a_2}+1)+\\
				&+\sum_{
					\begin{subarray}{c}
						a_1+a_2+a_3=k\\
						b_1+b_2+b_3=\ell
					\end{subarray}}
				(\mu_{1}^{a_1}+\mu_{2}^{a_1})
				(\mu_{1}+\mu_{2})^{a_2}
				2^{b_2}
				(\mu_{1}^{a_3}+\mu_{2}^{a_3})\\
				&\sum_{
					\begin{subarray}{c}
						|\boldsymbol\alpha|=a_1+b_1+1\\
						|\boldsymbol\beta|=a_2+b_2+1\\
						|\boldsymbol\delta|=a_3+b_3+1
					\end{subarray}}
				\zeta(\alpha_0,\cdots,\alpha_{a_1},
					\beta_0,\beta_1,\cdots,\beta_{a_2}+
					\gamma_0,\gamma_1,\cdots,\gamma_{a_3}+1).
		\end{split}
		\end{equation}

	Applying Ohno's relation (\cite{O}) for the index 
		$
	(\underbrace{1,\cdots,1}_{a_1},2,
				\underbrace{1,\cdots,1}_{a_2},2)
		$
	, we can see that
		\begin{align*}
			&\sum_{
					\begin{subarray}{c}
						b_1+b_2=\ell
					\end{subarray}}
				\sum_{
					\begin{subarray}{c}
						|\boldsymbol\alpha|=a_1+b_1+1\\
						|\boldsymbol\beta|=a_2+b_2+1
					\end{subarray}}
				\zeta(\alpha_0,\alpha_1,\cdots,\alpha_{a_1}+1,
					\beta_0,\beta_1,\cdots,\beta_{a_2}+1)\\
			=&\sum_{
					\begin{subarray}{c}
						b_1+b_2=\ell
					\end{subarray}}
				\zeta(a_2+b_2+2,a_1+b_1+2).
		\end{align*}
	Then the first term on the right hand side of (\ref{eq:after}) is same 
	as the the second term on the left hand side of (\ref{eq:after}).
	The second term on the right hand side of (\ref{eq:after}) can be written as follows.

	\begin{align*}
		\sum_{b_1+b_2+b_3=\ell}
				\sum_{
					\begin{subarray}{c}
						|\boldsymbol\alpha|=a_1+b_1+1\\
						|\boldsymbol\beta|=a_2+b_2+1\\
						|\boldsymbol\delta|=a_3+b_3+1
					\end{subarray}}
				2^{b_2}\zeta(\alpha_0,\cdots,\alpha_{a_1},
					\beta_0,\beta_1,\cdots,\beta_{a_2}+
					\gamma_0,\gamma_1,\cdots,\gamma_{a_3}+1)\\
				=\sum_{|\boldsymbol\alpha|=k+\ell+3}
					2^{\alpha_{a_1+1}+\cdots+\alpha_{a_1+a_2+1}-1-a_2}
					(1-2^{1-\alpha_{a_1+a_2+1}})
					\zeta(\alpha_0,\alpha_1,\cdots,\alpha_{k+1}+1).
	\end{align*}
	Then we have

	\begin{align*}
		&\sum_{|\boldsymbol\alpha|=k+\ell+3}
		\sum_{a_1+a_2+a_3=k}
		(\mu_{1}^{a_1}+\mu_{2}^{a_1})
		(\mu_{1}+\mu_{2})^{a_2}
		(\mu_{1}^{a_3}+\mu_{2}^{a_3})\times\\
		&\times2^{\alpha_{a_1+1}+\cdots+\alpha_{a_1+a_2+1}-1-a_2}
		(1-2^{1-\alpha_{a_1+a_2+1}})
		\zeta(\alpha_1,\alpha_1,\cdots,\alpha_{k+1}+1)\\
		=&(\ell+1)\zeta(k+\ell+4)\sum_{a_1+a_2=k}
			\mu_{1}^{a_1}\mu_{2}^{a_2}.
	\end{align*}
	
	In particular, when $\mu_1=-\mu_2=\mu\neq0$ and $k$ is even, the right hand side is
	\[
		\mu^{k}(\ell+1)\zeta(k+\ell+4),
	\]
	and the left hand side is
	\[
		\mu^{k}\left\{
		\sum_{|\boldsymbol\alpha|=k+\ell+3}
		\sum_{j=0}^{\frac{k}{2}}2^{\alpha_{2j+1}+1}
		\zeta(\alpha_{0},\alpha_{1},\cdots,\alpha_{k+1}+1)
		-2(k+2)\zeta(k+\ell+4)\right\}.
	\]
	Therefore we have
	\[
		\sum_{|\boldsymbol\alpha|=k+\ell+3}
		\sum_{j=1}^{\frac{k}{2}+1}2^{\alpha_{2j+1}+1}
		\zeta(\alpha_{1},\alpha_{2},\cdots,\alpha_{k+2}+1)
		=(2k+\ell+5)\zeta(k+\ell+4).
	\]

	For positive integers $m,n$ with $m>2n$, setting $k=2n-2$, $\ell=m-k-3$,
	we get \cite[Main Theorem]{ELO}.
	When $m=2n$, \cite[Main Theorem]{ELO} is derived from the duality theorem.
	Therefore we now conclude that \cite[Main Theorem]{ELO} can be deduced from our
	{\bf Theorem \ref{theorem:sub theorem}}.

	\begin{remark}
		More generally, by calculating integrals
			\begin{align*}
				\displaystyle &I_{k,\ell}(\mu_1,\cdots,\mu_n,\xi_1,\cdots,\xi_n)\\
		=&\frac{1}{k!\ell!}\mathop\int\limits_{
			\begin{subarray}{c}
						0<x_1<y_1<1\\
						\vdots\\
						0<x_n<y_n<1
			\end{subarray}}
	\left(\sum_{i=1}^n\mu_i\log\frac{1-x_i}{1-y_i}\right)^k
	\left(\sum_{j=1}^n\xi_j\log\frac{y_j}{x_j}\right)^{\ell}
	\prod_{h=1}^n\frac{dx_{h}dy_{h}}{(1-x_{h})y_{h}},
			\end{align*}
		we could get the same type of relations, but it seems difficult to write down the general form explicitly.
	\end{remark}

\end{document}